\documentclass[preprint,12pt]{elsarticle}

\usepackage{amssymb}
\usepackage{amsmath}
\usepackage{amsfonts}
\usepackage{amssymb}

\newcommand{\cit}{\mathbb{C}}

\newcommand{\modu}[1]{\left\vert#1\right\vert}
\newcommand{\nor}[1]{\left\Vert#1\right\Vert}
\newcommand{\h}{\mathcal H}
\newcommand{\bh}{\mathcal{ B(H)}}

\newcommand{\wa}{\mathit W(A)}
\newcommand{\waxb}{\overline{\mathit W(A_{x})}}

\newcommand{\wab}{\overline{\mathit W(A)}}

\newcommand{\set}[1]{\{#1\}}

\newtheorem{theorem}{Theorem}[section]

\newtheorem{lemma}[theorem]{Lemma}
\newtheorem{corollary}[theorem]{Corollary}

\newtheorem{definition}[theorem]{Definition}

\newtheorem{remark}[theorem]{Remark}

\journal{Applied Mathematics Letters}

\begin{document}

\begin{frontmatter}



\title{An application of normaloid operators}



\author[mymainaddress]{Abderrahim Baghdad}

\cortext[mycorrespondingauthor]{Abderrahim Baghdad}
\corref{mycorrespondingauthor}
\ead{bagabd66@gmail.com}

\author[mysecondaryaddress]{Mohamed Chraibi Kaadoud}
\ead{chraibik@uca.ac.ma}

\address[mymainaddress]{Department  of Mathematics FSSM, Cady Ayyad University \\ Marrakech, Morocco}
\address[mysecondaryaddress]{Department  of Mathematics FSSM, Cady Ayyad University \\ Marrakech, Morocco}



\begin{abstract}
Let $ x=(x_{n})_{n} $ be a bounded complex sequence and let
$ M_{x}= \displaystyle \sup_{n} \modu{x_{n}}$.
 By using a normaloid operator related to the sequence  $ x=(x_{n})_{n} $, we prove that
 $$  \sup_{\lambda \in \mathbb{C},~ \modu{\lambda}\leq M_{x}} \sup_{n} \modu{x_{n}+\lambda}=2M_{x}.$$  
\end{abstract}


\end{frontmatter}


\section{Introduction}
\indent{•}
Let $ x=(x_{n})_{n} $ be a bounded complex sequence and let
$ M_{x}=\displaystyle \sup_{n} \modu{x_{n}}.$
 For all integer $n$ and scalar $\lambda$, with $\modu{\lambda}\leq M_{x}$,  we have $\modu{x_{n}+\lambda}\leq \modu{x_{n}}+\modu{\lambda}\leq 2M_{x}$. It results that 
 \begin{equation}\label{inegalite} 
 \displaystyle \sup_{\lambda \in \cit,~ \modu{\lambda}\leq M_{x}} \sup_{n} \modu{x_{n}+\lambda}\leq2M_{x}.
 \end{equation} 
However, the equality 
\begin{equation}\label{egalite} 
\displaystyle \sup_{\lambda \in \cit,~ \modu{\lambda}\leq M_{x}} \sup_{n} \modu{x_{n}+\lambda}=2M_{x} 
\end{equation}
 is less obvious.\\
 \indent
  In this paper, we prove the equality \eqref{egalite} by using a normaloid operator related to the sequence  $ x=(x_{n})_{n} $. For this purpose, we need some notations and results from the literature.\\
\indent
Let $\bh$ be the algebra of all bounded linear operators acting on a complex Hilbert space ($\h,\langle .,. \rangle $). For $ A \in \bh$, the  numerical range of $A$ is denoted and defined by
$$ \wa = \set {\langle Ax, x\rangle  : x\in \h,  \   \nor{x} = 1}.$$
It is a celebrated result due to Toeplitz and Hausdorff that $ \wa $ is a nonempty bounded convex subset of the complex plane (not necessarily closed), and its supremum modulus, denoted and given by $$ w(A)= \sup \set {\vert\lambda\vert:~\lambda \in \wa},$$ 
 is called the  numerical radius of $ A $.
It is well known  that $ w(.) $ defines a norm on $\bh$, which is equivalent to the $C^{*}$-norm $\nor{.}$. In fact, the following inequalities are well-known:
$$ \frac{1}{2}\nor{A} \leq w(A)  \leq  \nor{A}, $$
for all $A\in \bh$, and an  operator  $ A\in \bh $ is said to be normaloid if the second inequality becomes an equality, i.e., $w(A)= \nor{A}$.\\
 \indent
For more material about the numerical radius and other information on the basic theory of algebraic numerical range and normaloid operators, we refer the reader to the books \cite{Bonsalla1, Bonsalla2, Gustafson, Halmos}. \\
 \indent
 From now on,  $\bh$ will denote the algebra of all bounded linear operators acting on a complex Hilbert space  $ \h $.  We shall denote the complex numbers by $\cit$.
\section{Main results}
Let us start with a definition and a lemma that we will use for giving two characterizations of a normaloid operator. These characterizations will allow us to prove the equality~\eqref{egalite}.
\begin{definition}
Let $\Delta$ be a nonempty bounded subset of $\cit$. The supremum modulus of $\Delta$ is denoted by $ \modu{\Delta}$ and given by
$$ \modu{\Delta} = \sup \set {\vert\lambda\vert:~\lambda \in \Delta}.$$
In particular, if  $A \in \bh$, then $\modu{\wa}=\modu{\wab}=w(A)$, where $\overline{\wa}$ is the closure of $\wa$.
\end{definition}

\begin{lemma}\label{ss}
Let $\set{0} \neq K$ be a  nonempty   compact subset of  $\cit$. Then there exists $s \in K$, $s\neq 0$ such that
$$\modu{ K+s}=\modu{K}+\modu{s}.  $$
\end{lemma} 

By a compactness argument, there exists  $s\in K$, $s\neq 0$ such that $ \modu{K}=\modu{s}$. So, $\modu{K+s}\leq \modu{K}+\modu{s}=2\modu{s}$. But $ 2s\in K+s$, then $2\modu{ s} \leq \modu{K+s} $ and therefore $\modu{ K+s}=\modu{K}+\modu{s}$.

\begin{theorem} \label{anormaloid} 
Let  $A \in \bh$. Then  the following  are equivalent statements
\begin{enumerate}
 \item  $A$ is normaloid,     

\item  $ \nor{A+\lambda}=\nor{A}+\modu{\lambda} $ for some $\lambda \in \cit^{*}$.

\end{enumerate} 
\end{theorem}

{\bf Proof.}\\
$(1)\Rightarrow (2)$. Assume that  $A$ is normaloid; that is $\modu{\wa}=w(A)=\nor{A}$. Since $\overline{\wa}$ is a compact subset of $\cit$, by Lemma~\ref{ss}, there exists $\lambda \in \wab$ ($\modu{\lambda}=w(A)$) satisfying
$$  \modu{\wab+\lambda} =\modu{\wab} +\modu{\lambda}=\nor{A}+\modu{\lambda}. $$
But $\wab+\lambda=\overline{W(A+\lambda)}$, then  $\modu{\wab+\lambda}=\modu{\overline{W(A+\lambda)}}=w(A+\lambda) \leq \nor{A+\lambda}$. It results that $ \nor{A}+\modu{\lambda} \leq \nor{A+\lambda} $ and so, $ \nor{A+\lambda}=\nor{A}+ \modu{\lambda} $. \\ 
$(2)\Rightarrow (1)$. Let $\lambda$ be a nonzero scalar such that $ \nor{A+\lambda}=\nor{A}+\modu{\lambda} $. By Barraa-Boumazgour~\cite[Theorem  2.1]  {BB}  we have 
$$ \modu{\lambda}\nor{A} \in \overline{W(\overline{\lambda}A)}=\overline{\lambda}~\overline{W(A)}.$$
 It results that  $ \nor{A} \leq w(A)$  and hence $ \nor{A}= w(A)$. This is exactly to say that $A$ is normaloid.

\begin{corollary}\label{coro} 
Let $A\in \bh$. Then $A$ is normaloid if and only if $$\sup_{\lambda\in \wab}\nor{A+\lambda}=2\nor{A}.$$
 \end{corollary}

{\bf Proof.}\\
Note first that, since $\modu{\lambda}\leq \nor{A}$ for any $\lambda \in \wab$, we always have
$$ \displaystyle\sup_{\lambda\in \wab}\nor{A+\lambda}\leq \sup_{\lambda\in \wab}(\nor{A}+\modu{\lambda})\leq 2\nor{A}.$$
If $A$ is normaloid, then  Theorem~\ref{anormaloid} states that there exists $\lambda \in \wab$, $\modu{\lambda}=w(A)=\nor{A}$ and $ \nor{A+\lambda}=\nor{A}+\modu{\lambda} =2\nor{A}$. Hence,
$$ \displaystyle\sup_{\lambda\in \wab}\nor{A+\lambda}=2\nor{A}.$$
Conversely,  since $\wab$ is compact, there exists $\mu \in \wab$  such that $$ \sup_{\lambda\in \wab}(\nor{A}+\modu{\lambda})=\nor{A}+\modu{\mu}.$$
If $A$ is not normaloid, then $\modu{\mu}<\nor{A}$ and we obtain 
$$ \sup_{\lambda\in \wab}\nor{A+\lambda}=\nor{A}+\modu{\mu}  <2\nor{A}.$$

\begin{remark}\label{att}   
If the operator $A$ is normaloid, the supermum $2\nor{A}$ is attained by  "some" scalars $\lambda \in \wab$ with $\modu{\lambda}=\nor{A}$. Then, we can write
\begin{equation}\label{req} 
\sup_{\lambda\in \wab}\nor{A+\lambda}=\sup_{\lambda\in \wab,~\modu{\lambda}\leq \nor{A}}\nor{A+\lambda}  =2\nor{A}.
\end{equation}
\end{remark}

\noindent
Now, we are ready to prove the following theorem wich  establishes the desired equality~\eqref{egalite}.
\begin{theorem}  
Let $ x=(x_{n})_{n} $ be a bounded complex sequence. Then,
 $$ \displaystyle \sup_{\lambda \in \cit,~ \modu{\lambda}\leq M_{x}} \sup_{n} \modu{x_{n}+\lambda}=2M_{x},$$  
 where $ M_{x}=\displaystyle \sup_{n} \modu{x_{n}}.$
\end{theorem}

{\bf Proof.}\\
Let $A_{x}$ be the diagonal operator defined by
$$A_{x}=
\begin{pmatrix}
 x_{1}   &  0    & \dots             &     &   \\
0       &x_{2}  &  &  &        & \\
\vdots &  & x_{3}  &  &  &    \\
 &  &  &   &  &   \\
 &        & &  & \ddots & \\
 &      &        &      &      &
\end{pmatrix}.$$
It is clear that $A_{x}$ is normaloid, indeed $w(A_{x})=\nor{A_{x}}=M_{x}$. Accoroding to the equality ~\eqref{req}, we get
$$\sup_{\lambda\in \waxb,~ \modu{\lambda}\leq \nor{A_{x}}}\nor{A_{x}+\lambda}=2 \nor{A_{x}}.$$  
Since $\nor{A_{x}+\lambda} =\displaystyle \sup_{n} \modu{x_{n}+\lambda}$  for all $\lambda \in \cit$ and $\nor{A_{x}}=M_{x}$, we obtain
$$   \displaystyle \sup_{\lambda \in \waxb,~ \modu{\lambda}\leq M_{x}} \sup_{n} \modu{x_{n}+\lambda}=2M_{x}.$$ 
It suffices to use the inequality~\eqref{inegalite} and the fact that $\waxb$ is a subset of $\cit$ to conclude that 
$$ \displaystyle \sup_{\lambda \in \cit,~ \modu{\lambda}\leq M_{x}} \sup_{n} \modu{x_{n}+\lambda}= 2M_{x}.$$











\section*{References}

\begin{thebibliography}{00}


\bibitem{BB}  M. Barraa,  M. Boumazgour, {\it   Inner  derivations and norm equality}, Proc. Amer. Math. Soc.,  Number 2 (2002), 471--476. 
\bibitem{Bonsalla1}  F. F. Bonsall and J. Duncan, {\it Numerical ranges of operators on normed spaces and of elements of normed
algebras}, London Mathematical Society Lecture Note Series 2 Cambridge University Press, London-New
York, (1971).
\bibitem{Bonsalla2} F. F. Bonsall and J. Duncan,  {\it Numerical ranges II}, London Mathematical Society Lecture Notes Series
10 Cambridge University Press, New York-London, (1973).
 

 \bibitem {Gustafson}   K. E. Gustafson, D. K. M.  Rao,   {\it  Numerical range: The Field of Values of Linear Operators and Matrices},    New York,  NY, USA, (1997).

\bibitem  {Halmos} P.R. Halmos, A Hilbert Space Problem Book, Van Nostrand, New York, 1967.


\end{thebibliography}


\end{document}